# Individual departure time decision considering departure scheduling utility

ZHANG Wen-yi(张文义)[1], GUAN Wei(关伟)[1], SUN Hui-jun(孙会君)[1], MAO Bao-hua(毛保华)[1]

1. MOE Key Laboratory for Urban Transportation Complex Systems Theory and Technology, Beijing Jiaotong University, Beijing 100044, China

**Abstract:** The scheduling utility plays a fundamental role in addressing the commuting travel behaviours. In this paper, a new scheduling utility, termed as DMRD-SU, was suggested based on some recent research findings in behavioural economics. DMRD-SU admitted the existence of positive arrival-caused utility. In addition, besides the travel-time-caused utility and arrival-caused utility, DMRD-SU firstly took the departure utility into account. The necessity of the departure utility in trip scheduling was analysed comprehensively, and the corresponding individual trip scheduling model was presented. Based on a simple network, an analytical example was executed to characterize DMRD-SU. It can be found from the analytical example that: 1) DMRD-SU can predict the accumulation departure behaviors at NDT, which explains the formation of daily serious short-peak-hours in reality, while MRD-SU cannot; 2) compared with MRD-SU, DMRD-SU predicts that people tend to depart later and its gross utility also decrease faster. Therefore, the departure utility should be considered to describe the traveler's scheduling behaviors better.

## 1 Introduction

Trip scheduling problem of commuter was of fundamental importance to study peak-period traffic congestion and analyze traffic control as well as the broader demand-side congestion relief measures, such as pricing and ride sharing incentives [1]. During the past several decades, this topic had attracted great research interests, spreading across the economical modeling [2-8], dynamic network analysis [9-12], trip scheduling under uncertainty [8, 13], and analyzing the value and reliability of travel time [14-17], etc. Among these studies, the economical modeling of scheduling utility (SU) usually played a fundamental role.

Based on the earlier work of Gaver [18] and Vickrey [19], Small [2] proposed the classical schedule delay concept which was later on extended to the uncertain case by Noland and Small [6]. According to the schedule delay concept, there exist a preferred arrival time (PAT), and arriving earlier or later than PAT would lead to penalty. The SU of schedule delay concept consists of linear travel time and arrival penalties (including early arrival penalty and late arrival penalty). De Palma et al. [5] changed PAT, a point, into an interval which was analogous to the indifference band proposed by Mahmassani and Chang [20]. Hendrickson and Plank [3] embedded risk attitudes into Small's SU, where travelers were assumed to be risk-averse and risk-prone to early arrival and late arrival, respectively. Besides risk attitudes, Li et al. [13] further took the monetary cost into account. However, Li's SU embedded the risk attitudes into the expected

**Foundation item:** Projects (71131001, 71271023) supported by National Natural Science Foundation of China; Project (2012CB725403) supported by National Basic Research Program of China; Project (2011AA110303) supported by National High Technology Research and Development Program of China; Project (2012JBZ005) supported by Fundamental Research Funds for the Central Universities, China
**Received date:** 2013-10-28; **Accepted date:** 2014-01-15
**Corresponding author:** GUAN Wei, Professor, PhD; Tel: +86-10-51687142; E-mail: weig@bjtu.edu.cn

rather than intrinsic outcomes. These studies have one in common, i.e., denying any kind of positive arrival utility in SUs. However, Senbil and Kitamura [7] and Jou et al. [8] empirically verified the existence of positive arrival utility in SU and suggested formulating SU by the multi-reference-point value function (MRP-VF) of prospect theory (see Refs. [21, 22] for the standard prospect theory). For more details on SU, refer to Refs. [7, 17], for example.

Prospect theory (see Refs. [21, 22] for detail) was firstly proposed by Kahneman and Tversky in 1979. Compared with the traditional behavioral economic theory, prospect theory relaxed the perfect rationality assumption made on people. Due to the excellent explain-ability on people's bounded rationality behaviors in reality, prospect theory had created significant impact in micro-economical theory and Kahneman and his colleagues were awarded the Nobel Prize in Economics in 2002 for their innovative contributions. Recently, Koszegi and Rabin's studies [23, 24] further reported that, besides the gain-loss (or relative) utility argued by prospect theory, people still concerned with the intrinsic consumption (or absolute) utility when assessing an alternative. Hence, besides the relative utilities caused by early and late arrivals, the travel-time-caused intrinsic utility should be also incorporated into scheduling utility. However, MRP-VF-based SUs did not consider the travel-time-caused (absolute) utility but only the arrival-caused (relative) utility. In addition, to the best of our knowledge, none of the existing SUs had explicitly considered the departure utility although it was implied. Knight [25] regarded the marginal utility of time spend at home as an essential motivation of departure time decision. Consequently, a sound SU should still takes the departure utility into account. For this reason, this study suggests a departure-utility-included multi-reference-dependent SU (DMRD-SU) to capture a traveler's perceived departure utility more precisely.

Along this line, the remaining context is organized as follows. In Section 2, the necessity of departure utility in DMRD-SU is illustrated; the formulation of DMRD-SU and the individual trip scheduling model are also presented in this section. Section 3 displays an analytical example is performed to explore the effects of departure utility in commuting departure scheduling behavior; Section 4 concludes this study and suggests some valuable works in the future.

## 2 Departure-utility-included scheduling utility and individual trip scheduling model

In this section, to simplify the presentation, the abbreviations used throughout this study are listed. A hypothesized trip case is introduced to explain the necessity of departure utility in scheduling utility. Then the departure-utility-included multi-reference-dependent SU (DMRD-SU) is presented. In the end, the individual trip scheduling model is given.

### 2.1 Abbreviations

| | |
|---:|---|
| SU | scheduling utility |
| MRP-VF | multi-reference-point value function |
| MRP-SU | MRP-VF-based SU |
| MRD-SU | (departure-utility-excluded) multi-reference-dependent SU |
| DMRD-SU | departure-utility-included multi-reference-dependent SU |
| PAT | preferred arrival time |
| NDT | normal departure time |
| PAE | preferred earliest-arrival time |
| PAL | preferred latest-arrival time |
| UT | travel-time-caused utility |
| UD | departure-caused utility |

UA    utility caused by arrival
UAE   utility caused by arrival earlier than PAT
UAL   utility caused by arrival later than PAT
EDT   earliest departure time
LDT   latest departure time

**2.2 A simple case study**

Below we introduce a hypothesized trip case (see Case 1) to identify the combination of a proper SU.

**Case 1.** A man goes to work from home to company during two days: in the first day he departs at 7: 30 and arrives at 8:00; in the second day he departs at 7: 20 and arrives at 8:00 because of an accident.

**Analysis:** When the travel-time-caused intrinsic disutility is ignored, we can infer SU for two work trips are the same since the arrival times are identical, causing the 10 minutes of difference in travel time consumption has no effect on the perceived disutility. Obviously, it is unreasonable. Thus, the travel time disutility cannot be excluded from SU. If the intrinsic travel time disutility is integrated into MRP-SU, we derive MRD-SU. However, judging by MRD-SU, SUs for two days differ from each other just due to their different consumed travel times. In other words, the deviation in the travel time disutility seizes the entire deviation regardless of the surrounded environments (e.g., at home, in office, in board, etc.) and trip objectives. It also violates the scientific founding that the value of duration varies with the environments, objectives, and some other factors [26]. Consequently, there should be additional utility related to the departure itself in daily trips. Actually, from departing at 7:30 but arriving at the same time with starting at 7:20, what a traveller achieves is not only to avoid the disutility caused by 10min of congestion, noise, populated air and non-profit time consumption, but also to enjoy some extra welfare. More precisely, since staying at home has much more freedom than in board, one can gain more positive utility by doing something that cannot be done (well) in board (e.g., 10min of benefit exercises, 10min more with family, etc.). Shortly, for the travels do not aim at the travels themselves, the thrifts on travel time mean not just the decreased time-loss but also some extra positive utility. In order to capture these psychologies, the departure utility should be also taken into account.

Base on the above analysis, we present the subsequent argument for constructing a sound scheduling utility.

**Argument 1.** A sound scheduling utility should comprise the travel-time-caused intrinsic utility, the arrival utility, as well as the departure utility.

For this, we suggest incorporating departure utility into MRD-SU, deriving the DMRD-SU. Until now, we can summarize the differences among SUs (see Table 1).

**Table 1** Compositions of diverse SUs

| SUs | Incorporated or not? | | | |
|---|---|---|---|---|
| | UD | UT | UA | Positive UA |
| Small [2]; Noland and Small [6]; Li et al. [13] | NO | YES | YES | NO |
| Mahmassani and Chang [4]; De Palma et al. [5] | NO | YES | YES | NO |
| Hendrickson and Plank [3] | NO | YES | YES | NO |
| MRP-SUs [7, 8] | NO | NO | YES | YES |
| MRD-SU in this paper | NO | YES | YES | YES |
| DMRD-SU in this paper | YES | YES | YES | YES |

**2.3 Formulation of DMRD-SU**

To formulate DMRD-SU, we introduce a pseudo reference-point called the preferred arrival time (PAT), and three real reference-points which are orderly named the normal departure time (NDT), preferred earliest-arrival time (PAE) and preferred latest-arrival time (PAL). Given the departure time $s$ and experienced travel time T, the present DMRD-SU is formulated as

$$GU = UD + UT + UAE + UAL \\ = \varphi(s;\text{NDT}) + u(\text{T}) + \rho(s,\text{T};\text{PAE},\text{PAT}) + \psi(s,\text{T};\text{PAT},\text{PAL}). \quad (1)$$

It can be known from Table 1 that MRP-SU comprises UAE and UAL, and MRD-SU comprises UT, UAE and UAL. Here, UD, UAE and UAL are gain-loss (or relative) utilities, while UT is the intrinsic utility. According to the reference-dependent theory [23, 24], the gain-loss utility is exactly the relative utility of prospect theory [21, 22], while the intrinsic utility is the outcome-based utility. Accordingly, UD, UT, UAE and UAL in Eq. (1) can be formulated orderly as follows.

$$\varphi = \begin{cases} -\kappa_\varphi^1 (\text{NDT}-s)^{\alpha_\varphi^1}, & \forall s \leq \text{NDT} \\ \kappa_\varphi^2 (s-\text{NDT})^{\alpha_\varphi^2}, & \forall s > \text{NDT} \end{cases} \quad (2)$$

$$u(\text{T}) = -\kappa_T \text{T} \quad (3)$$

$$\rho = \begin{cases} -\kappa_\rho^1 (\text{PAE}-s-\text{T})^{\alpha_\rho^1}, & \forall s+\text{T} \in (-\infty, \text{PAE}] \\ \kappa_\rho^2 (s+\text{T}-\text{PAE})^{\alpha_\rho^2}, & \forall s+\text{T} \in (\text{PAE}, \text{PAT}] \end{cases} \text{ and 0 otherwise} \quad (4)$$

$$\psi = \begin{cases} \kappa_\psi^1 (\text{PAL}-s-\text{T})^{\alpha_\psi^1}, & \forall s+\text{T} \in (\text{PAT}, \text{PAL}] \\ -\kappa_\psi^2 (s+\text{T}-\text{PAL})^{\alpha_\psi^2} + \Delta, & \forall s+\text{T} \in (\text{PAL}, +\infty) \end{cases} \text{ and 0 otherwise} \quad (5)$$

Here, $\Delta \leq 0$ is the disutility when the arrival time exceeds PAL, the parameter vector $\alpha = (\alpha_\varphi^1, \alpha_\varphi^2, \alpha_\rho^1, \alpha_\rho^2, \alpha_\psi^1, \alpha_\psi^2) \in (0,1]$ reflects that the traveler is risk-aversion toward gains and risk-prone toward loss, and $\kappa = (\kappa_\varphi^1, \kappa_\varphi^2, \kappa_T, \kappa_\rho^1, \kappa_\rho^2, \kappa_\psi^1, \kappa_\psi^2) > 0$ denotes the preference to each attribute. Moreover, $\left(\dfrac{\kappa_\varphi^1}{\kappa_\varphi^2}, \dfrac{\kappa_\rho^1}{\kappa_\rho^2}, \dfrac{\kappa_\psi^2}{\kappa_\psi^1}\right) \geq 1$ is restricted here to capture the traveler's loss-aversion behavior. Since the intrinsic utility is an outcome-based utility, UT has distinct appearance with the other gain-loss utilities. The parameters in Eqs. (2-5) can be calculated by the state-preference or revealed-preference survey data in practice. Fig. 1 displays the structural diagram of DMRD-SU.

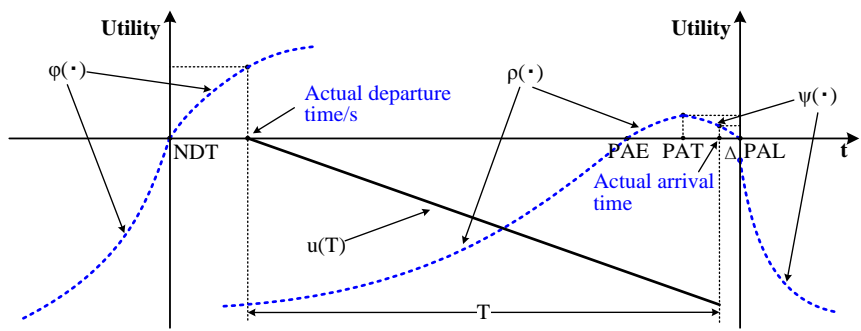

**Fig. 1** Schematic diagram of DMRD-SU

In Fig. 1, the thick dotted segments depict the additional utility created by departures and arrivals, and the thick sold segments describe the intrinsic consumption disutility of travel time. Note that there are two gain regions, i.e., the departures after NDT and the arrivals between PAE and PAL, whereas the others are loss ones. We hypothesize the traveler behaves risk-aversion towards the gains between PAE and PAL, and distinguish the two cases (i.e. the gains between PAE and PAT as well as PAT and PAL) through different degrees of loss-aversion (i.e., $\kappa_\rho^2 \leq \kappa_\psi^1$). Taking the UD segments away, the schematic diagram of MRD-SU is derived.

### 2.4 Individual trip scheduling model

From the perspective of an individual, everyone expect to obtain the most scheduling utility from his/her departure activity. Then the individual trip scheduling problem co-optimizing the route and departure time choices can be formulated as follows:

$$\max_{\theta,s} \sum_{k \in K} \theta_k \mathrm{GU}_k(s_k) \tag{6}$$

subject to

$$\sum_{k \in K} \theta_k = 1, \; \theta_k = 0 \text{ or } 1, \text{ and } s_k \in [\mathrm{EDT}, \mathrm{LDT}_k], \forall k \in K. \tag{7}$$

Here, $K$ is the feasible routes, $k$ is the route index, $\theta_k = 1$ if route $k$ is selected and 0 otherwise, EDT is the earliest departure time, $\mathrm{LDT}_k = \mathrm{PAL} - \mathrm{T}_k^{\inf}$ is the latest departure time, the first equation in Eq. (7) represents the traveler can only choose a single route for a departure action, and the second one is the feasible scheduling intervals. Consequently, the optimal route is the route with $\theta_k^* = 1$, and $s_k^*$ is the optimal departure time.

Studies [27, 28] had reported that, when encountering the congestion, the commuters tended to adjust the departure times far more readily than shift the routes. Accordingly, the previous trip scheduling model can be further simplified as follows:

$$\max_s \mathrm{GU}(s) \text{ subject to } s \in [\mathrm{EDT}, \mathrm{LDT}]. \tag{8}$$

## 3 Analytical example

In order to characterize the DMRD-SU and explore the impacts of departure utility on the individual departure time decision, a simple network with one route connecting a single OD pair (see Fig. 2) is applied to perform the subsequent analyses.

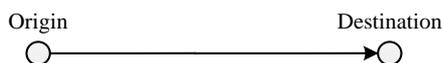

**Fig. 2** A simple network

This example is to explore the impacts of departure utility on individual scheduling behaviors. Hence, DMRD-SU and MRD-SU will be both applied. Moreover, since the assumption of diminishing sensitivity (see Ref. [21]) is too strict [23, 24] and to characterize the implication of loss aversion without diminishing sensitivity as a force on behavior, we set $\alpha = (1,1,1,1,1,1)$ here. Then we derive

$$\text{UD} = \begin{cases} -\kappa_\varphi^1 (\text{NDT} - s), \forall s \leq \text{NDT} \\ \kappa_\varphi^2 (s - \text{NDT}), \ \forall s > \text{NDT} \end{cases} \quad (9)$$

$$\text{UT} = -\kappa_T \text{T} \quad (10)$$

$$\text{UAE} = \begin{cases} -\kappa_\rho^1 (\text{PAE} - s - \text{T}), \forall s + \text{T} \in (-\infty, \text{PAE}] \\ \kappa_\rho^2 (s + \text{T} - \text{PAE}), \forall s + \text{T} \in (\text{PAE}, \text{PAT}] \end{cases} \text{ and 0 otherwise} \quad (11)$$

$$\text{UAL} = \begin{cases} \kappa_\psi^1 (\text{PAL} - s - \text{T}), \quad \forall s + \text{T} \in (\text{PAT}, \text{PAL}] \\ -\kappa_\psi^2 (s + \text{T} - \text{PAL}) + \Delta, \forall s + \text{T} \in (\text{PAL}, +\infty) \end{cases} \text{ and 0 otherwise.} \quad (12)$$

Based on above formulas, the DMRD-SU-based individual trip scheduling model can be formulated as

$$\max_s (\text{UD} + \text{UT} + \text{UAE} + \text{UAL}) \text{ subject to } s \in [\text{EDT}, \text{LDT}], \quad (13)$$

and the MRD-SU-based model can be formulated as

$$\max_s (\text{UT} + \text{UAE} + \text{UAL}) \text{ subject to } s \in [\text{EDT}, \text{LDT}]. \quad (14)$$

Based on above two models, we begin the analytical analyzes.

**Proposition 1.** Should a trip be scheduled, regardless of the congestion levels, we have $\kappa_\psi^2 \geq \kappa_\varphi^2$ for the DMRD-SU-based scheduling model.

**Proof.** Here we prove Proposition 1 by contradiction. Suppose $\kappa_\psi^2 < \kappa_\varphi^2$, we can conclude, when the departure time is later than $\text{PAL} - \text{T}$, the received GU will increase linearly with slope $\kappa_\varphi^2 - \kappa_\psi^2 > 0$ all the time, which implies the later departure the more one will gain. For this case, people will not depart, which contradicts with the fact that "trip was scheduled". Hence, Proposition 1 holds. □

Based on Proposition 1 and assume that $\kappa_\psi^2 > \kappa_\varphi^1 > \kappa_\psi^1 > \kappa_\varphi^2$ in MRD-SU and DMRD-SU in further, we derive another three propositions as follows.

**Proposition 2.** The following statements are valid for the DMRD-SU-based trip scheduling model:

$$s^* = \begin{cases} \text{PAT} - \text{T}, & \text{if } \text{T} < \text{PAT} - \text{NDT} \\ \text{NDT}, & \text{if } \text{T} \in [\text{PAT} - \text{NDT}, \text{PAL} - \text{NDT}] \\ \text{PAL} - \text{T}, & \text{if } \text{T} > \text{PAL} - \text{NDT} \end{cases}$$

**Proof.** From Eqs. (9-12), UD increases linearly with slopes orderly equaling to $\kappa_\varphi^1$ and $\kappa_\varphi^2$ for $s \leq \text{NDT}$ and $s > \text{NDT}$; UAE increases steadily when $s + \text{T} \leq \text{PAT}$; UAL decreases linearly with slopes orderly equaling to $-\kappa_\psi^1$ and $-\kappa_\psi^2$ for $\text{PAT} \leq s + \text{T} \leq \text{PAL}$ and $s + \text{T} > \text{PAL}$. Based on these properties, let us formally begin with the proofs.

(1) Given $\text{T} < \text{PAT} - \text{NDT}$: 1) when $s \leq \text{PAT} - \text{T}$, we have $s + \text{T} \leq \text{PAT}$, then the arrival utility, computed by UAE, increases steadily for $s + \text{T} \leq \text{PAT}$; 2) when $s > \text{PAT} - \text{T}$, we have $s + \text{T} > \text{PAT}$, the decrease rate of the arrival utility, computed by UAL, is faster than the increase rate of UD, which causes GU to reduce. Thus, $s^* = \text{PAT} - \text{T} > \text{NDT}$ if $\text{T} < \text{PAT} - \text{NDT}$.

(2) Given $\text{T} \in [\text{PAT} - \text{NDT}, \text{PAL} - \text{NDT}]$: 1) when $s < \text{NDT}$, we have $s + \text{T} < \text{PAL}$, no matter the arrival utility is computed by UAE or UAL, GU will increase since $\kappa_\varphi^1 > \kappa_\psi^1$; 2) when $s > \mathrm{NDT}$, we have $s + \mathrm{T} > \mathrm{PAT}$, the decrease rate of the arrival utility, computed by UAL, is faster than the increase rate of UD, which causes GU to reduce. Thus, $s^* = \mathrm{NDT}$ if $\mathrm{T} \in [\mathrm{PAT} - \mathrm{NDT}, \mathrm{PAL} - \mathrm{NDT}]$.

(3) Given $\mathrm{T} > \mathrm{PAL} - \mathrm{NDT}$: 1) when $s > \mathrm{PAL} - \mathrm{T}$, we have $s + \mathrm{T} > \mathrm{PAL}$, then the decrease rate of the arrival utility, computed by UAL, is faster than the increase rate of UD since $\kappa_\psi^2 > \kappa_\varphi^1 > \kappa_\varphi^2$, causing GU to reduce steadily; 2) however, when $s < \mathrm{PAL} - \mathrm{T}$, we have $s + \mathrm{T} < \mathrm{PAL}$, then no matter arrival utility is computed by UAE or UAL, GU will increase since $\kappa_\varphi^1 > \kappa_\psi^1$. Thus, $s^* = \mathrm{PAL} - \mathrm{T}$ if $\mathrm{T} > \mathrm{PAL} - \mathrm{NDT}$. □

Proposition 2 shows that, as the congestion level grows, the DMRD-SU-based model predicts earlier and earlier departure times, which agrees with the intuition. However, this relationship is not generally strict. There has also predicted a short-term departure convergence at NDT when the congestion is moderate, i.e., $\mathrm{T} \in [\mathrm{PAT} - \mathrm{NDT}, \mathrm{PAL} - \mathrm{NDT}]$. In this congestion interval, people commonly depart at the preferred departure time NDT, which may explain the formation causes of daily short-peak-hours. This also in return demonstrates the existence of NDT and departure utility hypothesized in this study.

**Proposition 3.** For the MRD-SU-based scheduling model, the optimal departure time will be strictly earlier and earlier as the congestion becomes worse and $s^* \equiv \mathrm{PAT} - \mathrm{T}$.

**Proof.** For the MRD-SU-based scheduling model, GU will be determined by the travel time disutility and arrival utility. Given the travel time T, GU will be solely determined by the arrival utility, and people will definitely arrived at the destination at PAT. Then we can deduce $s^* \equiv \mathrm{PAT} - \mathrm{T}$, which shows that the optimal departure time is strictly earlier and earlier as the congestion becomes worse. □

**Proposition 4.** The optimal departure times predicted by the DMRD-SU-based model are universally not earlier than that of the MRD-SU-based model for all congestion levels, and the optimal departure curves are plotted as follows (see Fig. 3):

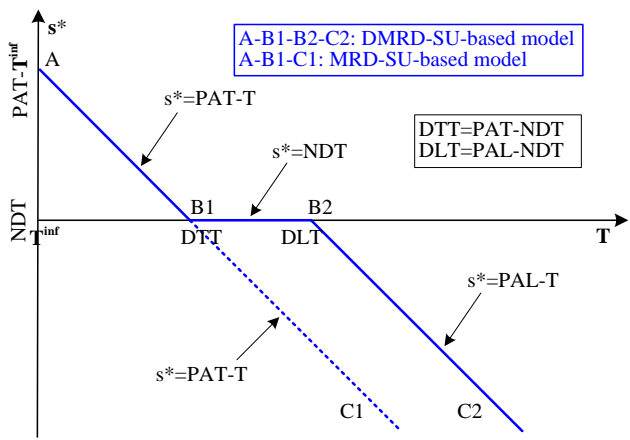

Fig. 3 The optimal departure times for two scheduling models

**Proof.** For the DMRD-SU-based scheduling model, it follows from Proposition 2 that

$$s^* = \begin{cases} \mathrm{PAT} - \mathrm{T}, & \text{if } \mathrm{T} < \mathrm{PAT} - \mathrm{NDT} \\ \mathrm{NDT}, & \text{if } \mathrm{T} \in [\mathrm{PAT} - \mathrm{NDT}, \mathrm{PAL} - \mathrm{NDT}] \\ \mathrm{PAL} - \mathrm{T}, & \text{if } \mathrm{T} > \mathrm{PAL} - \mathrm{NDT} \end{cases}$$

For the MRD-SU-based scheduling model, Proposition 3 concludes $s^* = \mathrm{PAT} - \mathrm{T}$. Then, Proposition 4 is proved and Fig. 3 can be derived. □

Fig. 3 indicates that, when the congestion is light (i.e., $T < PAT - NDT$), both models share the same optimal departure curves; however, from point-B1 on, after a stop between $PAT - NDT$ and $PAL - NDT$, DMRD-SU-based model predicts an universal departure lateness of $PAL - NDT$ than MRD-SU-based model. This should attribute to that UD can create extra positive utility for later departure.

Based on the above three propositions, we further plot the relationships (see Fig. 4) between the optimal gross disutility $GU^*$ and the congestion levels for both DMRD-SU-based and MRD-SU-based scheduling models. Note that the optimal gross disutility is assumed to be universally negative.

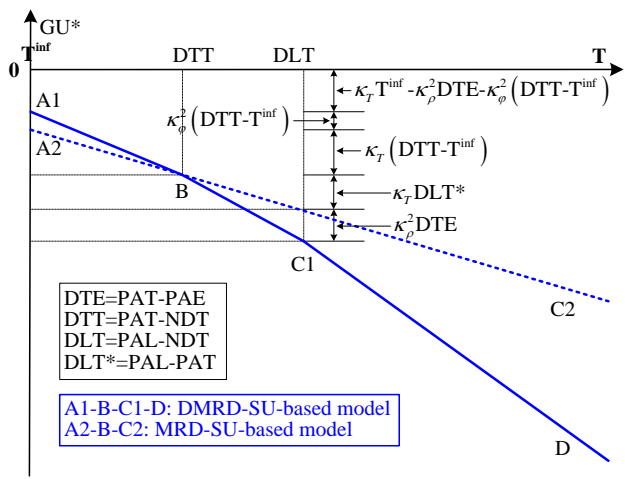

**Fig. 4** The optimal gross disutility for two scheduling models

Fig. 4 indicates that, the decrease rate of the optimal gross disutility $GU^*$ for DMRD-SU-based model is universally faster than that of MRD-SU-based model. This is due to, besides the travel time disutility, an UD-included traveler's loss still contains the (directly or indirectly) departure-caused disutility as the congestion becomes worse and worse. Also, when $T < PAT - NDT$, the UD-included $GU^*$ is larger than the UD-excluded one; when $T > PAT - NDT$, the UD-excluded $GU^*$ is larger. These are due to the positive UD created by the NDT-later departures and the negative UD caused by the NDT-earlier departures.

## 5 Conclusions and future researches

Based on the analysis from a hypothesized trip case, the authors argue that a sound scheduling utility should comprise the travel-time-caused intrinsic utility, the arrival utility, as well as the departure utility. For this, a new scheduling utility, termed as DMRD-SU, is suggested in this paper. Compare with the existing SUs, DMRD-SU firstly takes the departure utility into consideration. It is found from the simple analytic example on a simple network that:

1) DMRD-SU can predict the accumulation departure behaviors at NDT, which explains the formation of daily serious short-peak-hours in reality. In contrast, MRD-SU cannot.
2) Compared with MRD-SU, DMRD-SU predicts that people tend to depart later and its gross utility also decrease faster.

In conclusion, the departure utility should be considered to describe the traveler's scheduling behaviors better. Many further works are worthy of exploring based on the proposed DMRD-SU, for example,

1) The stochastic DMRD-SU and it-based individual (as well as aggregate) trip scheduling behaviors can be further studied.

  2) Empirical studies need to be executed to specify the behavioral parameters.
  3) Also, by adding monetary cost into the proposed DMRD-SU, the approach will be more realistic.